\documentclass[12pt]{article}

\usepackage[T2A]{fontenc}
\usepackage[cp1251]{inputenc}
\usepackage[english]{babel}
\usepackage{amsfonts}
\usepackage{eufrak}
\usepackage{amssymb}
\usepackage{amsmath}

\begin{document}

\newtheorem{theorem}{Theorem}[section]
\newtheorem{lemma}{Lemma}[section]
\newtheorem{corollary}{Corollary}[section]
\newtheorem{definition}{Definition}[section]
\newtheorem{proposition}{Proposition}[section]
\newtheorem{remark}{Remark}[section]

\begin{center}
	\textbf{\Large The Kagan characterization theorem on Banach spaces}
	
	\bigskip
	
	{\Large Margaryta Myronyuk}
	
	\bigskip
	
	\textit{B. Verkin Institute for Low Temperature Physics and Engineering\\ of the National Academy of
	Sciences of Ukraine, \\ Nauky Ave. 47, Kharkiv, Ukraine}
	
	\bigskip
	
	\textit{Bielefeld University, \\ Universit\"{a}tsstra\ss{}e 25, Bielefeld, Germany}
	
	\bigskip
	
	myronyuk@ilt.kharkov.ua
\end{center}

\begin{abstract}
	A. Kagan introduced classes of distributions $\mathcal{D}_{m,k}$ in $m$-dimensional space $\mathbb{R}^m$. He proved that if the joint distribution of $m$ linear forms of $n$ independent random variables belong to the class $\mathcal{D}_{m,m-1}$ then the random variables are Gaussian. If $m=2$ then the Kagan theorem implies the well-known Darmois-Skitovich theorem, where the Gaussian distribution is characterized by the independence of two linear forms of $n$ independent random variables. In the paper we describe Banach spaces where the analogue of the Kagan theorem is valid.
\end{abstract}



\bigskip

\section{Introduction}

In 1953 G. Darmois and V.P. Skitovich independently proved the following theorem which characterizes the Gaussian distribution on the real line.

\bigskip

\noindent\textbf{The Darmois-Skitovich theorem} (\cite{Dar}, \cite{Ski}, see also \cite[Sec. 3]{KaLiRa})
	\textit{Let $\xi_1,\dots,\xi_n$, $n \ge 2,$ be
	independent random variables. Let $a_j, b_j$ be nonzero numbers. If linear forms
	$L_1=a_1\xi_1+\dots+a_n\xi_n$ and $L_2=b_1\xi_1+\dots+b_n\xi_n$ are
	independent then $\xi_j$ are Gaussian random variables.} 

\bigskip

The generalizations of this characterization theorem to various algebraic structures were considered in numerous
works (see, e.g., \cite{Feldman2017} -- \cite{GhurO}, \cite{Myronyuk2006}-- \cite{Myronyuk2019}, \cite{Myronyuk2021},  for more references one can see \cite{Feldman2008book}), in which the coefficients of forms are topological automorphisms of the corresponding
algebraic structures.
In 1962, S.G. Ghurye and I. Olkin generalized the Darmois-Skitovich theorem to the case of random variables taking
values in $\mathbb{R}^n$, where the coefficients of forms are non-degenerate matrices (\cite{GhurO}). Later, A. Zinger simplified their
proof (see \cite[Sec. 3.2]{KaLiRa}). The Darmois-Skitovich theorem was also generalized to the case of random variables taking
values in a Hilbert space (see \cite{KannanKannappan}, \cite{Laha}). In 1985, W. Krakowiak generalized the Darmois-Skitovich theorem to the case of random variables taking
values in an
arbitrary separable Banach space $X$, where the coefficients of forms are invertible linear bounded
operators $A_j$ and $B_j$ (\cite{Krakowiak}). In 2008, M.Myronyuk gave another proof of the Darmois-Skitovich theorem on Banach spaces (\cite{Myr2008}). 
In this paper the Heyde theorem on Banach spaces was proved as well. Later in \cite{Myronyuk2021} the analogues of these theorems were proved for Q-independent random variables.

%Note that the proof of the Skitovich–Darmois theorem %and its generalizations reduces
%to the investigation of solutions of the following %equation in characteristic functionals:

Consider a random vector    $\boldsymbol{L}$= $(L_1, \dots, L_m)$. The distribution $\mu_{\boldsymbol{L}}$ of the random vector $\boldsymbol{L}$   belongs  to the class $\mathcal{D}_{m, k}$, $1\le k\le m$,  if the characteristic function $\hat\mu_{\boldsymbol{L}}(t_1, \dots,  t_m)$ admits the following factorization
$$
\hat\mu_{\boldsymbol{L}}(t_1, \dots,  t_m)={\mathbb{E}}[e^{{\mathrm i}\left\langle \boldsymbol{L}, (t_1, \dots,  t_m)\right\rangle }]=\prod_{{i_1, \dots, i_k}} R_{i_1, \dots i_k}(t_{i_1}, \dots t_{i_k}), \quad t_i\in \mathbb{R},
$$
where $R_{i_1, \dots, i_k}(t_{i_1}, \dots, t_{i_k})$ are continuous functions such that  $R_{i_1, \dots i_k}(0, \dots, 0)=1$, and all  indexes $(i_1, \dots, i_k)$ in this product  satisfy the conditions  $1\le i_1<\dots< i_k\le m$.
A. Kagan proved the following generalization of the Darmois-Skitovich theorem.

\bigskip

\noindent\textbf{The Kagan theorem} (\cite{Kagan1989})
\textit{Let $\xi_1,\dots,\xi_n$, $n \ge 2,$ be
	independent random variables. Let $A_{ji}$ be nonzero numbers. Put  $L_j=A_{j1}\xi_1+\cdots+A_{jn}\xi_n$, $j=1, 2, \dots, m$.  If
	$\boldsymbol{L}$= $(L_1, \dots, L_m) \in \mathcal{D}_{m, m-1}$ then all $\xi_j$ are Gaussian random variables.} 

\bigskip

\noindent The class $\mathcal{D}_{2, 1}$ contains all distributions of 
$\mathbb{R}^2$ with independent components. Therefore the Darmois-Skitovich is a partial case of the Kagan theorem. Note that the Kagan theorem was generalized for locally compact Abelian groups (\cite{Fe-Kagan-2020}). In this paper we prove  analogues of the Kagan theorem for Banach spaces.

\section{Notation and definitions}

We use standard facts of the theory of probability distributions in
a Banach space (see, e.g., \cite{Vah-Tar-Chob}). Let $X$ be a
separable real Banach space and let $X^*$ be its dual space.
Denote by ${\cal L}_{inv}(X)$ the set all linear bounded invertible
operators in $X$. Denote by $A^*$ the adjoint operator to a linear
bounded operator $A$. Denote by $I$ the identity operator. Denote
by $\left\langle x,f\right\rangle $ the value of a functional $f\in X^*$ at an element $x\in
X$. Denote by $\mathcal{M}^1(X)$ the convolution semigroup of
probability distributions on $X$. For $\mu \in \mathcal{M}^1(X)$ we
denote by $\overline{\mu}$ the distribution defined by the formula
$\overline{\mu}(E)=\mu(-E)$ for any Borel set $E\subset X$. Let
$\xi$ be a random variable with values in $X$ and with a
distribution $\mu$. The characteristic functional of $\mu$ is
defined by the formula

$$ \widehat\mu(f)=\mathbf{E}[e^{\mathrm{i}\left\langle \xi,f\right\rangle }]=\int\limits_X e^{\mathrm{i}\left\langle x,f\right\rangle }
d\mu(x), \quad f\in X^*.$$

\noindent Note that
$\widehat{\overline{\mu}}(f)={\overline{\widehat\mu}}(f)$. It is
known that $\widehat\mu(f)$ is positive definite and continuous in
the norm topology (\cite[Ch.4, \S2]{Vah-Tar-Chob}).

In the paper, unless otherwise specified, we work in a norm
topology.

Let $\psi(f)$ be a complex-valued function on $X^*$, and let $h$ be
an arbitrary element of $X^*$. We denote by $\Delta_h$ an operator
of the finite difference
$$\Delta_h \psi(f)=\psi(f+h)-\psi(f).$$
A function $\psi(f)$ on $X^*$ is called a polynomial if
$$\Delta_{h}^{n+1}\psi(f)=0$$
for some $n$ and for all $f,h \in X^*$.  

%The minimal $n$ for which
%this equality holds is called the degree of the polynomial
%$\psi(f)$.

%\begin{definition} \label{0.12} A distribution $\mu \in %\mathcal{M}^1(\mathbb{R})$ is called a Gaussian (or normal) distribution on %$\mathbb{R}$ if either $\mu$ is a degenerated distribution or $\mu$
%is absolutely continuous with respect to the Lebesgue measure in $\mathbb{R}$ and %has the density
%$$\rho(t)={1\over \sqrt{2\pi\sigma}}e^{-{(t-m)^2\over
%2\sigma}}, \quad t\in\mathbb{R}, \quad m\in\mathbb{R}, \quad\sigma>0.$$
%The characteristic function of a Gaussian distribution has the form 
%$$\hat\mu(t)= e^{ {i mt-{1\over 2}\sigma t^2 }}, \quad t\in\mathbb{R},
%\quad m\in\mathbb{R}, \quad \sigma\geq 0.$$
%\end{definition}

\begin{definition} \label{0.2} A random variable $\xi$ is called Gaussian if
 for any $f\in X^*$ the real-valued random variable
$\left\langle \xi,f\right\rangle $ is Gaussian. In other words $\mu \in \mathcal{M}^1(X)$ is called Gaussian if its one-dimensional image $\mu_f=\mu\circ f^{-1}$ is a Gaussian distribution for all  $f\in X^*$ and hence there exist real numbers $m_f$ and $\sigma_{f}\geq 0$ such that $\hat\mu_f(t)= e^{ {\mathrm{i} m_f t-{1\over 2}\sigma_f t^2 }}$, $t\in\mathbb{R}$.
\end{definition}

%Note that the following statement is true in a Banach space (see (\cite[Ch.4, %\S2.4]{Vah-Tar-Chob})):

%\begin{proposition} 
%Let $X$ be a Banach space. The characteristic functional of a Gaussian %distribution $\mu\in
%\mathcal{M}^1(X)$ has the form
%\begin{equation}\label{0.3}
%    \hat\mu(f)= e^{ i<m,f>-{1\over 2}<R f,f>},\ f\in X^*,
%\end{equation}
%where $m\in X$ and $R: X^*\longrightarrow X$ is a symmetric nonnegative operator. %The element $m$ is called the mean
%value of the distribution $\mu$ and $R$ is called its covariation operator.

%Conversely, if $\mu\in \mathcal{M}^1(X)$ has a characteristic functional of the %form $(\ref {0.3})$, where $m\in X$ is a certain
%element and $R: X^*\longrightarrow X$ is a certain symmetric nonnegative %operator, then $\mu$ is a Gaussian distribution
%in $X$ with mean value $m$ and covariation operator $R$.

%\end{proposition}

%Taking into account that $\hat\mu(tf)=\hat\mu_f(t)$ for all $f\in
%X^*$ and $t\in \mathbb{R}$, we obtain

%\begin{equation}\label{0.4}
%    <m,f>=m_f, \quad <R f,f>=\sigma_f.
%\end{equation}

The following proposition follows immediately from Definition \ref{0.2}.

\begin{proposition} \label{0.5} A distribution $\mu$ on $X$ is
Gaussian if and only if for any $f\in X^*$ the characteristic
function $\widehat\mu_f(t)$ of a random variable $\left\langle \xi,f\right\rangle $ is a
characteristic function of a Gaussian distribution on $\mathbb{R}$.
\end{proposition}

Proposition \ref{0.5} allows easily to obtain analogues of the well-known classical theorems for a real separable Banach space such as the Kac-Bernstein theorem (\cite{Bernstein}, \cite{Kac}) on characterization of a Gaussian distribution by the independence of the sum of independent random variables and their difference, the Cramer theorem on the decomposition
of a Gaussian distribution (\cite[\S4.1]{Ramachandran}), and the Marcinkiewicz theorem (\cite[\S3.13]{Ramachandran}). The proofs of these
theorems reduces to the description of solutions of the corresponding functional equations that admit restrictions to each one-dimensional subspace. 

For the convenience of references we formulate the analogues of the Cramer theorem and the Marcinkiewicz theorem for Banach spaces as propositions.

\begin{proposition}\label{Cr}  Let $\xi$ be a Gaussian random
	variable with values in a real separable Banach space $X$. Let
	$\xi=\xi_1+\xi_2$, where $\xi_1, \xi_2$ are independent random
	variables with values in  $X$. Then $\xi_1, \xi_2$ are Gaussian
	random variables.
\end{proposition}

\begin{proposition}\label{Mar}  Let $X$ be a real
	separable Banach space, $\mu\in \mathcal{M}^1(X)$. Let $\psi(f)$ be
	a polynomial on $X^*$. If
	\begin{equation}\label{1.3}
	\widehat\mu(f)=e^{\psi(f)}
	\end{equation}
	\noindent in a certain neighborhood of zero, then $\mu$ is a
	Gaussian distribution on $X$.
\end{proposition}

Note that the anologue of the Darmois-Skitovich theorem for a Banach space can not be proved in a such simple way.

\subsection{The main results}

Let $X$ be a separable real Banach space. Consider a random vector    $\boldsymbol{\zeta}$= $(\zeta_1, \dots, \zeta_m)$ in values in $X^m$. Following A.M. Kagan (\cite{Kagan1989}) we say that the distribution of the random vector $\boldsymbol{\zeta}$   belongs  to the class $\mathcal{D}_{m, k}$, $1\le k\le m$,  if the characteristic functional $\hat\mu_{\boldsymbol{\zeta}}(f_1, \dots,  f_m)$ admits the following factorization
$$
\hat\mu_{\boldsymbol{\zeta}}(f_1, \dots,  f_m)={  \mathbb{E}}[e^{{\mathrm i}\left\langle \boldsymbol{\zeta}, (f_1, \dots,  f_m)\right\rangle }]=\prod_{{i_1, \dots, i_k}} R_{i_1, \dots i_k}(f_{i_1}, \dots f_{i_k}), \quad f_i\in X^*,
$$
where $R_{i_1, \dots i_k}(f_{i_1}, \dots, f_{i_k})$ are continuous functions such that  $R_{i_1, \dots i_k}(0, \dots, 0)=1$, and in the product  all  indexes $(i_1, \dots, i_k)$ satisfy the conditions  $1\le i_1<\dots< i_k\le m$.

We formulate an analogue of
this characterization for independent random variables with
values in a Banach space.

\begin{theorem}\label{KaganTh}
Let $X$ be a real reflexive separable Banach space.
Let $A_{pj}$, $p=\overline{1,m}$, $i=\overline{1,n}$, be linear continuous operators of $X$ such that $Ker A_{pj}=\{0\}$. Put

\begin{equation}\label{Kag2}
	G_i=\{ (A_{1i}x, ..., A_{mi}x) \in X^m: x\in X \},\quad i=\overline{1,n}.
\end{equation}

\noindent Suppose that the condition 

\begin{equation}\label{Kag1}
\overline{G_i} \cap \overline{G_j} = \{0\}, \quad i,j=\overline{1,n},
i\neq j.
\end{equation}

\noindent holds. Let $\xi_i$, $i=1, 2, \dots, n$,  be independent random variables with values in $X$. Consider the linear forms
$L_j=A_{j1}\xi_1+\cdots+A_{jn}\xi_n$, $j=1, 2, \dots, m$. If the random vector   $\mathbf{L}$=$(L_1, \dots, L_m)$ 
$\in\mathcal{D}_{m, m-1}$ then all $\xi_i$ are Gaussian random variables.
\end{theorem}	

Theorem \ref{KaganTh} is a spesial case of Theorem \ref{QKaganTh} which is formulated and is proved further for $Q$-independent random variables.
	
%\section{Q-analogues of the Kagan theorem}

\medskip

Let $\xi_1, \dots, \xi_n$ be random variables with values in
separable Banach space $X$. Following A. Kagan and G. Sz\'{e}kely
(\cite{KS}), we say that random variables $\xi_1, \dots, \xi_n$ are
$Q$-independent  if the characteristic functional of the
distribution of the vector $(\xi_1, \dots, \xi_n)$ can be
represented in the form
\begin{equation}\label{i0}
    \widehat\mu_{(\xi_1, \dots, \xi_n)}(f_1, \dots, f_n)={\bf E}[e^{i\left\langle (\xi_1, \dots, \xi_n),(f_1, \dots, f_n)\right\rangle }]
    =$$$$=\left(\prod_{j=1}^n\widehat\mu_{\xi_j}(f_j)\right)\exp\{q(f_1,
\dots, f_n)\}, \quad f_j\in X^*,
\end{equation}

\noindent where $q(f_1, \dots, f_n)$ is a continuous polynomial on $(X^*)^n$
such that $q(0, \dots, 0)=0$.

The paper \cite{KS} has led to new researches of characterization problems on locally compact Abelian groups (\cite{Feldman2017}, \cite{Myronyuk2019}, \cite{Myr2020}) and Banach spaces (\cite{Myronyuk2021}).

It is obvious that if random variables $\xi_1, \dots, \xi_n$ are
independent then they are $Q$-independent (in this case the function $q\equiv 0$).
Therefore if a statement is valid for Q-independent random variables then it is valid for independent random variables (for example Theorem \ref{QKaganTh} implies Theorem \ref{KaganTh}). 

%\subsection{The Cramer theorems (Q-version)}

The following $Q$-analogue of the Proposition \ref{Cr} (the Cramer theorem in a Banach space)
follows directly from Definition \ref{0.2} and the $Q$-analogue of the
Cramer theorem on a real line (see \cite{KS}).

\begin{proposition}\label{KrQ}
	Let $\xi$ be a Gaussian random
	variable with values in a real separable Banach space $X$. Let
	$\xi=\xi_1+\xi_2$, where $\xi_1, \xi_2$ are $Q$-independent random
	variables with values in  $X$. Then $\xi_1, \xi_2$ are Gaussian
	random variables.
\end{proposition}

The following theorem is the $Q$-analogue of Theorem \ref{KaganTh}.

\begin{theorem}\label{QKaganTh}
	Let $X$ be a real reflexive separable Banach space.
	Let $A_{pj}$, $p=\overline{1,m}$, $i=\overline{1,n}$, be linear continuous operators of $X$ such that $Ker A_{pj}=\{0\}$. Let $G_i$ be the sets defined by $(\ref{Kag2})$.
		%Put
	%\begin{equation*}
	%G_i=\{ (A_{1i}x, ..., A_{mi}x) \in X^m: x\in X \},\quad i=\overline{1,n}.
	%\end{equation*}
	Suppose that the condition $(\ref{Kag1})$ 
%	\begin{equation}\label{QKag1}
%	\overline{G_i} \cap \overline{G_j} = \{0\}, \quad i,j=\overline{1,n},
%	i\neq j.
%	\end{equation}
	holds. Let $\xi_i$, $i=1, 2, \dots, n$,  be Q-independent random variables with values in $X$. Consider the linear forms
	$L_j=A_{j1}\xi_1+\cdots+A_{jn}\xi_n$, $j=1, 2, \dots, m$. If the random vector   $\mathbf{L}$=$(L_1, \dots, L_m)$ 
	$\in\mathcal{D}_{m, m-1}$ then all $\xi_i$ are Gaussian random variables.
\end{theorem}

\textbf{Proof.}	Consider the continuous linear functions
$T_i: X \rightarrow X^m$ of the form $$T_i(x)=(A_{1i}x, ..., A_{mi}x)\in G_i$$ and the factor mappings $$P_i: X^m \rightarrow X^m/
\overline{G_i}.$$ We will construct the compositions $$F_{ij}=P_i T_j$$
and will check that $Ker F_{ij}=\{0\}$ for $i\neq j$. Let $x\in Ker
F_{ij}$. Then $F_{ij}(x)=P_i T_j(x)=P_i (A_{1j}x, ..., A_{mj}x) \in
Ker P_i$ that is $(A_{1j}x, ..., A_{mj}x) \in \overline{G_i}$. It follows from condition (\ref{Kag1}) that $x=0$ for $i\neq j$. Thus $Ker
F_{ij}=\{0\}$ for $i\neq j$. Since $X$ is a reflexive space, it follows from this that

\begin{equation}\label{QKagt7}
\overline{F_{ij}^*(X^*)}=X^*
\end{equation}

\noindent in the norm topology.

Note that $F_{ij}^*=T_j^* P_i^* $. Put $C_{pi}=A_{pi}^*$,
$p=\overline{1, m}$, $i=\overline{1, n}$. It is easy to verify that the mapping $T_i^*: (X^*)^m \rightarrow X^*$ has the form

$$ T_i^*(f_1,...,f_m)= C_{1i}f_{1}+ \cdots +C_{mi}f_{m}, \quad
i=\overline{1,n}.  $$
Put $$(\overline{G_i})^{\perp}=\{ (f_1,...,f_m)\in (X^*)^m:
\langle x_1,f_1 \rangle+ \ldots +\langle x_m,f_m \rangle=0 \quad
\forall (x_1,..., x_m)\in \overline{G_i} \}.$$ It is clear that $$(\overline{G_i})^{\perp}=\{ (f_1,...,f_m)\in
(X^*)^m: \langle A_{1i}x,f_1 \rangle+ \ldots +\langle A_{mi}x,f_m
\rangle=0 \quad \forall x\in X\}=$$$$= \{ (f_1,...,f_m)\in (X^*)^m:
\langle x, C_{1i}f_{1}+ \cdots +C_{mi}f_{m}\rangle=0 \quad \forall
x\in X\}=Ker T_i^* .$$ Thus $P_i^*$ is a natural embedding of
$Ker T_i^*$ into $(X^*)^m$.

Let $\mu_i$ be a distribution of the random variable $\xi_i$.
We consider the characteristic functional  $\hat\mu_{\mathbf{L}}(f_1,
\dots, f_m)$ of the random vector $\mathbf{L}$. Since the random variables $\xi_i$ are Q-independent, the functional $\hat\mu_{{\mathbf L}}(f_1,
\dots, f_m)$ has the form
$$\hat\mu_{{\mathbf L}}(f_1, \dots, f_m)={\bf E}[e^{{{\mathrm i}\left\langle {\mathbf L},(f_1, \dots,
f_m)\right\rangle }}]={\bf E}[e^{{{\mathrm i}\left\langle (L_1,\dots,L_m),(f_1, \dots,
f_m)\right\rangle }}]=$$
$$={\bf E}\bigg[ e^{{\mathrm i}\sum_{j=1}^m\left\langle A_{j1}\xi_1+\cdots+A_{jn}\xi_n,
	f_j\right\rangle }\bigg]={\bf E}\bigg[e^{{\mathrm i}\sum_{i=1}^n\left\langle \xi_i,
C_{1i}f_1+\dots+C_{mi}f_m\right\rangle }\bigg]=$$
\begin{equation}\label{QKagt01}
=\prod_{i=1}^n\hat\mu_i(C_{1i}f_1+\dots+C_{mi}f_m) \exp\{q(f_1,
\dots, f_m)\}, \quad f_j\in X^*.
\end{equation}

\noindent where $q(f_1, \dots, f_m)$ is a continuous polynomial on $(X^*)^n$
such that $q(0, \dots, 0)=0$.,

Since ${\mathbf L}\in\mathcal{D}_{m,
	m-1}$, we have

\begin{equation}\label{QKagt02}
	\hat\mu_{{\mathbf L}}(f_1, \dots, f_m)=\prod_{j=1}^{m}
	R_j(f_1,   \dots, f_{j-1},f_{j+1},\dots, f_m), \quad f_j\in X^*,
\end{equation}

\noindent where $R_{i_1, \dots i_k}(f_{i_1}, \dots, f_{i_k})$ are continuous functions such that  $R_{i_1, \dots i_k}(0, \dots, 0)=1$, and all  indexes $(i_1, \dots, i_k)$ satisfy the conditions  $1\le i_1<\dots< i_k\le m$.

It follows from (\ref{QKagt01}) and (\ref{QKagt02}) that
$$\prod_{i=1}^n\hat\mu_i(C_{1i}f_1+\dots+C_{mi}f_m)\exp\{q(f_1,
\dots, f_m)\}=$$
\begin{equation}\label{QKagt1}
=\prod_{j=1}^{m}
R_j(f_1,   \dots, f_{j-1},f_{j+1},\dots, f_m), \quad f_j\in X^*.
\end{equation}

Put $\nu_i=\mu_i*\bar\mu_i$. Then
$\hat\nu_i(f)=|\mu_i(f)|^2\geq 0$ and the characteristic functionals $\hat\nu_i(f)$ also satisfy (\ref{QKagt1}), in which all factors in the left and right parts of the equation are non-negative. If we prove that all $\nu_i$ are Gaussian distributions then Proposition \ref{KrQ} implies that all $\mu_i$ are Gaussian distributions. Thus, we may assume from the beginning that all factors in the left and right parts of the equation (\ref{QKagt1}) are non-negative. 

Since $\hat\mu_i(0)=1$, it follows from the continuity of the characteristic functionals $\hat\mu_i(f)$ in the norm topology that there exists a ball
$U_R$ around zero such that all $\hat\mu_i(f)> 0$ for $f\in U_R$.
We choose in $U_R$ a ball $V_r$, where the radius $r={R\over m \cdot \max\|
	C_{pi}\|}$ ($p=\overline{1, m}$, $i=\overline{1, n}$). Then

\begin{equation}\label{t2}
C_{1i}(V_r)+\dots+C_{mi}(V_r)\subset U_R, \quad i=\overline{1,n}.
\end{equation}

Put $\psi_i=\log \hat\mu_i$ in $V_r$,
$i=\overline{1,n}$, $s_j=\log R_j$, $j=1, 2, \dots, m$. It follows from (\ref{QKagt1}) that the functions $\psi_i$ and $s_j$ satisfy equation 

$$\sum_{i=1}^n \psi_i( C_{1i}f_1+\cdots+C_{mi}f_m  ) + q(f_1,
\dots, f_m) =$$
\begin{equation}\label{QKagt2}
=
\sum_{j=1}^m s_j(f_1, \ldots, f_{j-1}, f_{j+1}, \ldots, f_m),
\quad f_j\in V_r,
\end{equation}

\noindent where $s_j$ are arbitrary functions.

We use the finite difference method to solve equation (\ref{QKagt2}). Let $g_{nj}$,
$j=\overline{1,m}$, be arbitrary elements of $V_r$ such that
$$ C_{1n}g_{n1}+ \cdots +C_{mn}g_{nm}=0. $$
Note that the existence of such elements follows from (\ref{QKagt7}) and the form of $F_{ij}^*$. Substitute $f_j+g_{nj}$ for $f_j$ for all $j$ in equation (\ref{QKagt2}). We obtain

$$\sum_{i=1}^{n-1} \Delta_{h_{in}} \psi_i( C_{1i}f_1+\cdots+C_{mi}f_m  ) +\Delta_{\overline{{g}}_{n}} q(f_1,
\dots, f_m) =$$
\begin{equation}\label{QKagt3}
=
\sum_{j=1}^m \Delta_{\overline{{v}}_{jn}} s_j(f_1, \ldots, f_{j-1}, f_{j+1}, \ldots, f_m),
\quad f_j\in V_r,
\end{equation}

\noindent where $h_{in}= C_{1i}g_{n1}+ \cdots +C_{mi}g_{nm}$,
$i=\overline{1,n-1}$, $\overline{{v}}_{jn}= (g_{n1}, \ldots,
g_{n,j-1}, g_{n,j+1}, \ldots, g_{nm})$, $j=\overline{1,m}$, $\overline{{g}}_{n}=(g_{n1},...,g_{nm})$. The left-hand side of equation (\ref{QKagt3}) no longer contains the function $\psi_n$.

Let $g_{n-1,j}$, $j=\overline{1,m}$, be arbitrary elements of
$V_r$ such that

$$ C_{1,n-1}g_{n-1,1}+ \cdots +C_{m,n-1}g_{n-1,m}=0. $$
Reasoning similarly we exclude the function $\psi_{n-1}$ from the left-hand side of equation (\ref{QKagt3}). We obtain

$$\sum_{i=1}^{n-2} \Delta_{h_{in}} \Delta_{h_{i,n-1}}  \psi_i( C_{1i}f_1+\cdots+C_{mi}f_m  ) +\Delta_{\overline{{g}}_{n}} \Delta_{\overline{{g}}_{n-1}} q(f_1,
\dots, f_m) =$$
\begin{equation}\label{QKagt3.1}
=
\sum_{j=1}^m \Delta_{\overline{{v}}_{jn}} \Delta_{\overline{{v}}_{j,n-1}} s_j(f_1, \ldots, f_{j-1}, f_{j+1}, \ldots, f_m),
\quad f_j\in V_r,
\end{equation}

\noindent where $h_{i,n-1}= C_{1i}g_{n-1,1}+ \cdots +C_{mi}g_{n-1,m}$,
$i=\overline{1,n-2}$, $\overline{{v}}_{j,n-1}= (g_{n-1,1}, \ldots,
g_{n-1,j-1}, g_{n-1,j+1}, \ldots, g_{n-1,m})$, $j=\overline{1,m}$, $\overline{{g}}_{n-1}=(g_{n-1,1},...,g_{n-1,m})$.

By excluding the functions $\psi_{j}$ from the left-hand side of equation (\ref{QKagt2}), after $n-1$ steps we come to the equation of the form

$$\Delta_{h_{12}} ...\Delta_{h_{1n}} \psi_1( C_{11}f_1+\cdots+C_{m1}f_m  ) +
\Delta_{\overline{{g}}_{n}}...\Delta_{\overline{{g}}_{2}} q(f_1,\dots, f_m)=$$

\begin{equation}\label{QKagt4}
=
\sum_{j=1}^m r_j(f_1, \ldots, f_{j-1}, f_{j+1}, \ldots, f_m),
\quad f_j\in V_r,
\end{equation}

\noindent where

\begin{equation}\label{QKagt5}
h_{1j}= C_{11}g_{j1}+ \cdots +C_{m1}g_{jm}, \quad
j=\overline{2,n},
\end{equation}

\noindent $\overline{{g}}_{j}=(g_{j1},...,g_{jm})$ and $r_j$ are arbitrary functions. Note that we chose elements $g_{ij}$  at every step in such a manner that the equlity

\begin{equation}\label{QKagt6}
C_{1i}g_{i1}+ \cdots +C_{mi}g_{im}=0, \quad
i=\overline{2,n},
\end{equation}

\noindent is fulfilled. 

Let $k_{m}$ be an arbitrary element of $X^*$. Substitute $f_m+k_m$ for $f_m$ in (\ref{QKagt4}) and subtract equation (\ref{QKagt4}) from the resulting equation. 
We get 

\begin{eqnarray}\label{QKagt8}
\Delta_{C_{mi}k_m} \Delta_{h_{12}} ...\Delta_{h_{1n}} \psi_1( C_{11}f_1+\cdots+C_{m1}f_m  ) +
\Delta_{\overline{{k}}_{m}} \Delta_{\overline{{g}}_{n}}...\Delta_{\overline{{g}}_{2}} q(f_1,\dots, f_m)=\nonumber\\=
\sum_{j=1}^{m-1} r_j(f_1, \ldots, f_{j-1}, f_{j+1}, \ldots, f_m),
\quad f_j\in V_r,
\end{eqnarray}

\noindent where $\Delta_{\overline{{k}}_{m}}=(0,..,0,k_m)$. Note that equation (\ref{QKagt8}) does not contain the function $r_m$.

By repeating this operation, we consistently exclude all functions $r_j$ from the right-hand side of the resulting equations. After $m$ steps we get

\begin{eqnarray}\label{QKagl2.2}
\Delta_{C_{11}k_1}\dots \Delta_{C_{m1}k_m}\Delta_{h_{12}}
...\Delta_{h_{1n}} \psi_1( C_{11}f_1+\cdots+C_{m1}f_m  )+ \nonumber\\+
\Delta_{\overline{{k}}_{1}}...\Delta_{\overline{{k}}_{m}} \Delta_{\overline{{g}}_{n}}...\Delta_{\overline{{g}}_{2}} q(f_1,\dots, f_m)=0, \quad
f_j\in V_r.
\end{eqnarray}

Since $q(f_1, \dots, f_m)$ is a polynomial, we have

\begin{equation}\label{QKagt2.2}
\Delta^{l+1}_{(s_1,...,s_m)}q(f_1, \dots, f_m)=0,
\quad f_j\in X^*,
\end{equation}

\noindent for some $l$ and arbitrary elements $s_1,...,s_m \in X^*$.

Taking into account (\ref{QKagt2.2}), it follows from (\ref{QKagl2.2}) that

\begin{eqnarray}\label{QKagl2.2.1}
\Delta^{l+1}_{C_{11}s_1+\cdots+C_{m1}s_m} 
\Delta_{C_{11}k_1}\dots \Delta_{C_{m1}k_m}\Delta_{h_{12}}
...\Delta_{h_{1n}} \psi_1( C_{11}f_1+\cdots+C_{m1}f_m  )=0, \quad
f_j\in V_r.
\end{eqnarray}

Since $X$ is a reflexive space and $Ker A_{pj}=\{0\}$, the subgroups ${C_{j1}(X^*)}$ are dence in $X^*$ in the norm topology. Therefore we can put
$C_{11}k_1=\cdots=C_{m1}k_m=g$, where $g\in V_r$. Since elements $h_{ij}$
are defined by (\ref{QKagt5}) and belong to $F_{ij}^*(X^*)$, it follows from (\ref{QKagt7}) that we can put
$h_{12}=\cdots=h_{1n}=g$. Since $C_{11}s_1+\cdots+C_{m1}s_m$ belong to $F_{ij}^*(X^*)$, it follows from (\ref{QKagt7}) that we can put $C_{11}s_1+\cdots+C_{m1}s_m=g$. Thus, it follows from (\ref{QKagl2.2.1}) that 
the function $\psi_1(f)$ satisfies the equation

\begin{equation}\label{l2.3}
\Delta^{l+m+n}_{ g } \psi_1(f)=0, \quad   f, g\in V_r,
\end{equation}

\noindent  i.e. $\psi_1(f)$ is a polynomial in a neighborhood of zero of $X^*$. In the same way we can prove that all functions $\psi_j( f)$, $j=\overline{2,n}$, are polynomials in a neighborhood of zero of $X^*$. It follows from Proposition \ref{Mar}
that all $\mu_i$ are Gaussian distributions. $\square$

\medskip

%\subsubsection{Comments}

It is known that a characteristic functional is continuous in the norm topology and is sequentially continuous in the topology of pointwise convergence  (\cite[P.4, \S2]{Vah-Tar-Chob}). We use only the continuity in the norm topology to prove Theorem \ref{QKaganTh}. An essential fact, which was used in the proof, is the following assertion:
the equality $Ker A=\{0\}$ for some bounded
operator $A$ implies that $\overline{Im A^*}=X^*$. This assertion is not valid for arbitrary Banach spaces, but it is true for reflexive ones.

Suppose that $X$ is not a reflexive space. For which maximally wide class of Banach spaces does Theorem \ref{QKaganTh} remain true?  

Let us also assume at first that the random variables $\xi_i$,
$i=1, 2, \dots, n$, have the distributions $\mu_i$ such that their characteristic functions do not vanish. Then equation (\ref{QKagl2.2.1}) is fulfilled on $X^*$. Putting
$f_1=f, f_2=\cdots=f_m=0$, $s_1=s, s_2=\cdots=s_m=0$ in (\ref{QKagl2.2.1}), we get

\begin{eqnarray}\label{QKagl2.2.2}
\Delta^{l+1}_{C_{11}s} 
\Delta_{C_{11}g_1}\dots \Delta_{C_{m1}g_m}\Delta_{h_{12}}
...\Delta_{h_{1n}} \psi_1( C_{11}f_1)=0, \quad
f\in X^*.
\end{eqnarray}

\noindent Elements $C_{11}f$, $C_{11}s$, $C_{11}g_1, \ldots,C_{m1}g_m$, $h_{12}, \ldots, h_{1n}$ belong to images of the operators conjugate to operators with zero kernel. Since $X$ is not a reflexive space, we can not obtain from (\ref{QKagl2.2.2}) the equation of type (\ref{l2.3})
on $X^*$.

Now we will use the fact that a characteristic functional is  sequentially continuous in the topology of pointwise convergence. It means that if a sequence of elements $\{f_n\}$ of $X^*$ converges pointwise to $f\in X^*$ ($\{f_n(x)\rightarrow f(x)\}$) for any $x\in X$ then $\hat\mu(f_n)\rightarrow \hat\mu(f)$. We need the following definition:

\begin{definition}
	A Banach space $X$ is said to be quasireflexive if $X$ has finite codimension in its second dual $X^{**}$.
\end{definition}

%\medskip

%\textbf{Задача. } Нехай $X$ --- дійсний сепарабельний банахів
%простір, $A$ --- неперервний оператор $X$ такий, що $Ker A=\{0\}$.
%Коли вірно (для яких просторів $X$), що для будь-якого елементу
%$f\in X^*\backslash Im(A^*)(X^*)$ існує послідовність $\{f_n\}$
%елементів образу $Im(A^*)$, які поточково збігаються до $f$?

%\medskip

We need the following proposition.

\begin{proposition}\label{quasireflexive}
	Let $X$ be a quasireflexive separable Banach space. Let $A$ be a linear bounded operator of $X$ such that $Ker A=\{0\}$. Let $f$ be an arbitrary element of $ X^*\backslash Im(A^*)(X^*)$. Then there exists a sequence of elements $\{f_n\}$
	of the image $Im(A^*)$ such that a sequence of elements $\{f_n\}$ converges pointwise to $f$.
\end{proposition}
As V. Kadec informed the author, Proposition \ref{quasireflexive} holds (comments of V. Kadec see in Appendix). Note that readers  can also find useful information about quasireflexive spaces and additional references in Appendix.

Using Proposition \ref{quasireflexive} we obtain from (\ref{QKagl2.2.2}) the following equation:

\begin{equation}\label{k2}
\Delta^{l+m+n}_{ g } \psi_1(f)=0, \quad   f, g\in X^*.
\end{equation}

Now we can formulate the following theorem.

\begin{theorem}\label{QKaganTh2}
	Let $X$ be a real quasireflexive separable Banach space.
	Let $A_{pj}$, $p=\overline{1,m}$, $i=\overline{1,n}$, be linear continuous operators of $X$ such that $Ker A_{pj}=\{0\}$. Let $G_i$ be the sets defined by $(\ref{Kag2})$.
		Suppose that the condition $(\ref{Kag1})$ 
		holds.
		 Let $\xi_i$, $i=1, 2, \dots, n$,  be Q-independent random variables with values in $X$ and with non vanishing characteristic functionals. Consider the linear forms
	$L_j=A_{j1}\xi_1+\cdots+A_{jn}\xi_n$, $j=1, 2, \dots, m$. If the random vector   $\mathbf L$=$(L_1, \dots, L_m)$ belongs to the class
	$\mathcal{D}_{m, m-1}$ then all $\xi_i$ are Gaussian random variables.
\end{theorem}

Note that we cannot omit the condition that the characteristic functionals of the distributions $\mu_i$ do not vanish because otherwise we cannot claim that the equation (\ref{k2}) is fulfilled in some ball.

Theorem \ref{QKaganTh2} implies the following theorem.

\begin{theorem}\label{KaganTh2}
	Let $X$ be a real quasireflexive separable Banach space.
	Let $A_{pj}$, $p=\overline{1,m}$, $i=\overline{1,n}$, be linear continuous operators of $X$ such that $Ker A_{pj}=\{0\}$. Let $G_i$ be the sets defined by $(\ref{Kag2})$.
	Suppose that the condition $(\ref{Kag1})$ 
	holds.
	Let $\xi_i$, $i=1, 2, \dots, n$,  be independent random variables with values in $X$ and with non vanishing characteristic functionals. Consider the linear forms
	$L_j=A_{j1}\xi_1+\cdots+A_{jn}\xi_n$, $j=1, 2, \dots, m$. If the random vector   $\mathbf L$=$(L_1, \dots, L_m)$ belongs to the class
	$\mathcal{D}_{m, m-1}$ then all $\xi_i$ are Gaussian random variables.
\end{theorem}

\section{Appendix}

The content of this section belongs to V. Kadec and is published with his consent.

\subsection{Density and sequential density in weak-star topology}

Let $X$ and $Y$ be normed spaces over the field $\mathbb K$ ($\mathbb K=\mathbb R$ or $\mathbb K=\mathbb C$). We denote by $\mathcal{L}(X,Y)$ the space of all bounded linear operators from $X$ into $Y$. By $B_X$ and $S_X$ we denote the closed unit ball and the unit sphere of $X$ respectively. 
%Remaining notation is also standard and follows \cite{Kad} (mainly Chapter 17 of %that book). 

\begin{definition} 
	Let $X$ be a Banach space. A set $F \subset {X}^*$ is called \emph{total over} $X$, if for any $y\in X \setminus \{0\}$ there exists an $f \in F$ such that $f(y) \ne 0$.
\end{definition}
Evidently, $F$ is total over $X$ if and only if its linear span ${\rm Lin\ } F$ is total, and if and only if the norm-closure $\overline{\rm Lin\ } F$ of the linear span is total. The bipolar theorem applied to the dual pair $(X^*, X)$ implies that $F \subset {X}^*$ is total over $X$ if and only if ${\rm Lin\ } F$ is $w^*$-dense in $X^*$.

\begin{definition}
	Let $X$ be a Banach space and $\theta \in (0,1]$. A set $F \subset X^*$ is said to be $\theta$-\emph{norming over} $X$ if
	\begin{equation*}
	\sup_{f \in F\backslash \{0\}} \dfrac{|f(x)|}{\|f\|} \ge \theta \|x\| %\tag{1}
	\end{equation*}
	for all $x \in X$. The set $F \subset X^*$ is said to be \emph {norming} if there exists a $\theta \in (0,1] $ such that $F$ is $\theta$-norming over~$X$.
\end{definition}

The following well-known result can be deduced from the bipolar theorem, applied to the dual pair $(X^*, X)$.

\begin{proposition}[{\cite[Section 17.2.4, Exercise 2]{Kad}}] \label{prop-norming}
	Let $X$ be a Banach space, and $\theta \in (0,1]$. A set $F \subset S_{X^*}$ is $\theta$-norming for $X$ if and only if the $w^*$-closure of the absolute convex hull of $F$ contains $\theta {B}_{X^*}$.
\end{proposition}

\begin{proposition}[{\cite[Section 17.2.4, Corollary 4]{Kad}}] \label{prop-w*metrizable}
	Let $X$ be a separable Banach space. Then the 
	$w^*$-topology $\sigma(X^*, X)$ is metrizable on the bounded subsets of the space $X^*$. 
\end{proposition}
%\begin{definition} 
%	A Banach space $X$ is said to be \emph{quasireflexive} if $X$ has finite %codimension in its second dual $X^{**}$.
%\end{definition}

The very first example of quasireflexive nonreflexive space was constructed by R.C.James in \cite{James}. Since then, quasireflexive spaces form an important part of Banach space theory.

The following result by W.J.~Davis and J.~Lindenstrauss (\cite{Dav-Lin}) is highly non-trivial. Two years later it was rediscovered by A.~Plichko (\cite{Plich}). The proof in both papers develops ideas by Yuri Petunin (\cite{Pet}). One can find a detailed exposition with all necessary preliminaries in the book \cite{Pet-Pl}.

\begin{proposition} \label{prop-tot-not norming} 
	A Banach space $X$ has a total nonnorming subspace in its dual if and only if $X$ is not quasireflexive. 
\end{proposition}

\begin{definition} \label{def-seq-dense}
	Let $X$ be a Banach space. A set $F \subset {X}^*$ is called \emph{w*-sequentially dense} in $X^*$, if for every $g \in X^*$ there exists a sequence $f_n \in F$ such that $\lim_{n \to \infty}f_n(x) = g(x)$ for all $x \in X$.
\end{definition}

\begin{lemma} \label{lemma-sequen-norming} 
	Let $X$ be a separable Banach space, $F \subset X^*$ be a linear subspace. Then the following assertions are equivalent: $(i)$ $F$ is w*-sequentially dense in $X^*$, and $(ii)$ $F$ is norming. 
\end{lemma}

\textbf{Proof.}
	$(i) \Rightarrow (ii)$. Assume $F$ is w*-sequentially dense in $X^*$. Denote $U$ the w*-closure of the unit ball $B_F$ of $F$. Thanks to Banach-Steinhaus theorem \cite[Section 10.4.2, Theorem 1]{Kad}, the sequence $f_n$ from Definition \ref{def-seq-dense} is bounded, hence it belongs to a set of the form $n B_F$. Consequently, every $g \in X^*$ belongs to some set of the form $n U$, that is
	\begin{equation} \label{eq-nU} 
	X^* = \bigcup_{n \in \mathbb{N}} n U.
	\end{equation}
	The application of Baire theorem gives that $U$ contains some ball $\theta B_{X^*}$. Since $U$ is the $w^*$-closure of the absolute convex hull of $S_F$, Proposition \ref{prop-norming} says that $S_F$ is norming, so $F$ is norming as well.
	
	$(ii) \Rightarrow (i)$. Evidently, a linear subspace $F$ is norming if and only if its unit sphere is norming, so in notations of the previous implication Proposition \ref{prop-norming} says that $U$ contains some ball $\theta B_{X^*}$. So, \eqref{eq-nU} holds true. Consequently, every $g \in X^*$ belongs to the w*-closure of some ball $n B_{F}$. But, thanks to metrizability (Proposition \ref{prop-w*metrizable}), this means that $g$ is a w*-limit of some sequence from $n B_{F} \subset F$. $\square$

\medskip

After this, Proposition \ref{prop-tot-not norming} reformulates as follows:
\begin{theorem} \label{theor-sequen-norming2} 
	For a separable Banach space $X$ the following assertions are equivalent: $(i)$ every $w^*$-dense subspace of $X^*$ is $w^*$-sequentially dense, and $(ii)$ $X$ is quasireflexive. 
\end{theorem}

\subsection{Applications to the range of adjoint operator}

Recall that for arbitrary Banach spaces $X, Y$ and for arbitrary $T \in \mathcal{L}(X,Y)$ the $w^*$-density of $T^*(Y^*)$ in $X^*$ is equivalent to the injectivity of $T$ (\cite[Section 17.1.3, Theorem 5]{Kad}). This and Theorem \ref{theor-sequen-norming2} give the following corollary:

\begin{theorem} \label{theor-dual-op-range1} 
	Let $X$ be a separable quasireflexive Banach space, $Y$ be a Banach space and $T \in \mathcal{L}(X,Y)$ be injective. Then $T^*(Y^*)$ is $w^*$-sequentially dense in $X^*$.
\end{theorem}

It is not clear for us whether for every norm-closed total subspace $F \subset X^*$ there exist a Banach space $Y$ and an injective operator $T \in \mathcal{L}(X,Y)$ such that $T^*(Y^*) \subset F$. In the case of positive answer one would have the following inverse to Theorem \ref{theor-dual-op-range1} result:

\medskip

\textbf{Hypothesis}. Let $X$ be a separable Banach space that is not quasireflexive. Then there exist a Banach space $Y$ and an injective operator $T \in \mathcal{L}(X,Y)$ such that $T^*(Y^*)$ is not $w^*$-sequentially dense in $X^*$.

\medskip

We can demonstrate a partial result. At first, a technical tool.

\begin{lemma} \label{theor-dual-op-range2} 
	Let $X$, $Y$ be infinite-dimensional Banach spaces, and let $F \subset X^*$ be a $w^*$-separable closed linear subspace. Then there exists a $T \in \mathcal{L}(X,Y)$ such that $T^*(Y^*) \subset F$ and $T^*(Y^*)$ is $w^*$-dense in $F$.
\end{lemma}

\textbf{Proof.}
	Fix a $w^*$-dense in $S_F$ countable set $\{f_n\}_{n \in {\mathbb{N}}} \subset S_F$ and a normalized basic sequence $\{e_n\}_{n \in {\mathbb{N}}} \subset S_Y$ (the existence of the latter was demonstrated by S.~Masur, see a detailed exposition in \cite[Theorem 6.3.3]{KadKad}). Denote $\{e_n^*\}_{n \in {\mathbb{N}}} \subset Y^*$ there corresponding coordinate functionals extended to the whole $Y$ by the Hahn-Banach theorem. Define $T \in \mathcal{L}(X,Y)$ by formula
	$$
	Tx = \sum_{n=1}^\infty 2^{-n}f_n(x) e_n.
	$$
	Then, the adjoint operators $T^* \colon Y^* \to X^*$ acts on every $g \in Y^*$ as
	$$
	T^*g = \sum_{n=1}^\infty 2^{-n} g(e_n) f_n.
	$$
	Evidently, $T^*(Y^*) \subset \overline{{\rm Lin\ }}\{f_n\}_{n \in \mathbb{N}} \subset F$. On the other hand, for every $k \in \mathbb{N}$ we have $T^*e_k^* = 2^{-k} f_k$, consequently $T^*(Y^*)$ contains the $w^*$-dense in $F$ set ${\rm Lin\ } \{f_n\}_{n \in \mathbb{N}}$. $\square$

\begin{corollary} \label{corcor} 
	$X$, $Y$ be infinite-dimensional Banach spaces, $X$ and $X^*$ are both separable, and $X$ is not quasireflexive. Then there exist an injective operator $T \in \mathcal{L}(X,Y)$ such that $T^*(Y^*)$ is not $w^*$-sequentially dense in $X^*$.
\end{corollary}

\textbf{Proof.}
	Indeed, Theorem \ref{theor-sequen-norming2} gives the existence of a $w^*$-dense norm-closed linear subspace $F \subset X^*$ which is not $w^*$-sequentially dense. By separability of $X^*$ the subspace $F$ is separable, and consequently $F$ is $w^*$-separable. Applying Lemma \ref{theor-dual-op-range2} we obtain a $T \in \mathcal{L}(X,Y)$ such that $T^*(Y^*) \subset F$ and $T^*(Y^*)$ is $w^*$-dense in $F$. Then $T^*(Y^*)$ is $w^*$-dense in $X^*$, so $T$ is injective. But the inclusion $T^*(Y^*) \subset F$ implies that $T^*(Y^*)$ is not $w^*$-sequentially dense in $X^*$.
$\square$

Remark, that for separable $X$ the dual space $X^*$ is $w^*$-separable (\cite[Section 17.2.4, Corollary 2]{Kad}). Unfortunately, this does not guaranty the $w^*$-separability of subspaces of $X^*$. This is the reason, why in the Corollary \ref{corcor} we make the assumption of separability of $X^*$ in norm topology. We don't know if in the case of separable non-quasiraflexive $X$ the corresponding $w^*$-dense but not $w^*$-sequentially dense $F \subset X^*$ can be always selected to be $w^*$-separable. We cannot extract this additional property neither from the demonstration by Davis-Lindenstrauss (\cite{Dav-Lin}), nor from Plichko's demonstration, but such a possibility does not contradict our intuition. 

\section*{Acknowledgements}

The author would like to thank the Volkswagen Foundation (VolkswagenStiftung), the Bielefeld University and Prof. Dr. Friedrich G\"{o}tze for the support and warm reception.

\section*{Funding}

This research was supported by VolkswagenStiftung - Az. 9C108.

\end{document}